\documentclass{amsart}%
\usepackage{amssymb}
\usepackage{amsmath}
\usepackage{graphicx}
\usepackage{amsfonts}
\usepackage{color}
\usepackage{graphicx}%
\providecommand{\U}[1]{\protect\rule{.1in}{.1in}}
\newtheorem{theorem}{Theorem}
\theoremstyle{plain}

\numberwithin{equation}{section}

\begin{document}
\title{Edge Tessellations and Stamp Folding Puzzles}
\author{Matthew Kirby}
\address{Hershey High School\\
Hershey, PA 17033}
\email{maskirby@gmail.com}
\author{Ron Umble}
\address{Department of Mathematics\\
Millersville University of Pennsylvania\\
Millersville, PA 17551}
\email{ron.umble@millersville.edu}
\date{August 30, 2011}

\begin{abstract}
An edge tessellation is a tiling of the plane generated by reflecting a
polygon in its edges. We prove that a polygon generating an edge tessellation
is one the following eight types: a rectangle; an equilateral, 60-right,
isosceles right, or 120-isosceles triangle; a 120-rhombus; a 60-90-120 kite;
or a regular hexagon. A stamp folding puzzle is a paper folding problem
constrained to the perforations on a sheet of postage stamps. We establish the
following conjecture due to G. Frederickson: \textquotedblleft Although
triangular stamps have come in a variety of different triangular shapes, only
three shapes seem suitable for [stamp] folding puzzles: equilateral, isosceles
right triangles, and 60-right triangles.\textquotedblright\ 

\end{abstract}
\maketitle

Which polygons generate a tiling of the plane when reflected in their edges?
The complete answer, discovered by Millersville University students Andrew
Hall, Joshua York, and the first author in the spring of 2009, and we present
it here as a theorem:

\begin{theorem}
\label{one}A polygon generating a tiling of the plane when reflected in its
edges is one of the following eight types: a rectangle; an equilateral,
$60$-right, isosceles right, or $120$-isosceles triangle; a $120$-rhombus; a
$60$-$90$-$120$ kite; or a regular hexagon.
\end{theorem}

A \emph{tessellation }(or\emph{ tiling})\textit{ }of the plane is a collection
of plane figures that fills the plane with no overlaps and no gaps. An
\emph{edge tessellation }is generated by reflecting a polygon in its edges.
The eight edge tessellations, which are pictured in Figures 1 and 2, are the
most symmetric examples of \emph{Laves tilings}, as one can see from the
complete list in \cite{Gr-Sh}.
\begin{figure}[ptb]%
\centering
\includegraphics[
height=2.9395in,
width=4.9727in
]%
{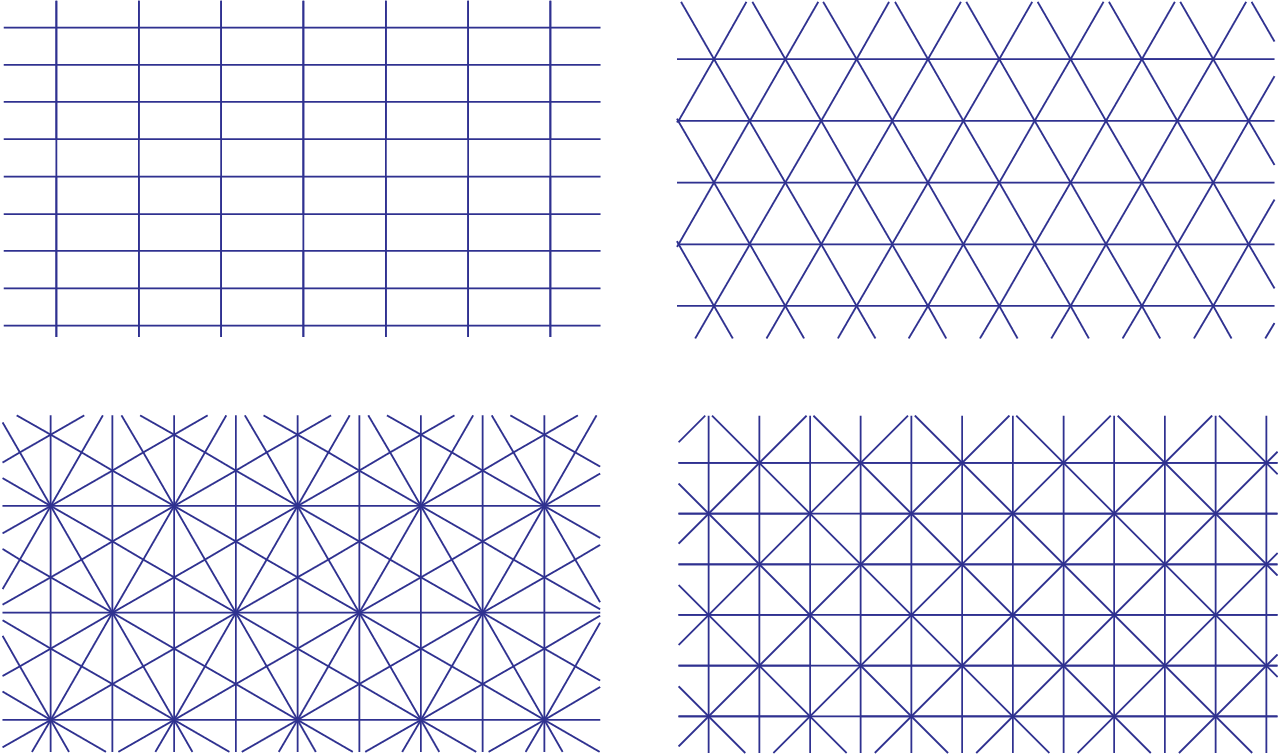}%
\caption{Figure 1. Edge tessellations generated by non-obtuse polygons.}%
\end{figure}
\begin{figure}[ptb]%
\centering
\includegraphics[
height=2.9326in,
width=4.9666in
]%
{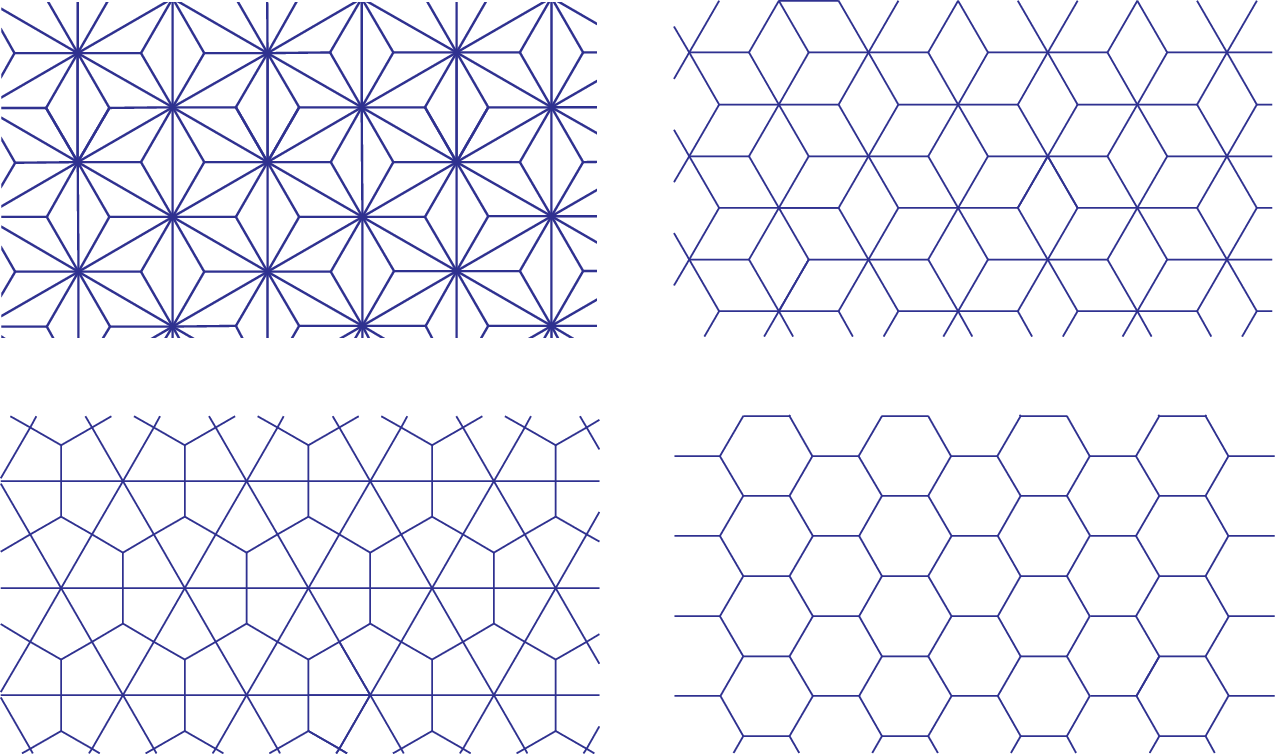}%
\caption{Figure 2. Edge tessellations generated by obtuse polygons.}%
\end{figure}

Edge tessellations provide the setting for \emph{stamp folding puzzles}, which
are paper folding problems constrained to the perforations on a sheet of
postage stamps. The sheet must embed in a edge tessellation, may have any
shape, and may be bounded or unbounded, with bounding edges along perforations
as in Figure 3. The object of a stamp folding puzzle is to create some
specified configuration by folding the sheet along its perforations without
creasing the stamps. Tucks, which slip one subpacket of folded stamps between
the leaves of another, are allowed. Indeed, a tuck is required to solve the
following delightful problem posed by G. Frederickson on page 144 of his book
\textquotedblleft Piano-Hinged Dissections:\ Time to Fold!\textquotedblright%
\textit{ \cite{Fr}}:\ \emph{Consider the block of sixteen isosceles right
triangular stamps pictured in Figure 3. Fold the block into a packet
sixteen-deep so that the stamps are arranged in the order} 4 1 16 6 5 15 14 8
7 13 11 12 2 3 9 10.%
\begin{figure}[ptb]%
\centering
\includegraphics[
height=1.574in,
width=1.5705in
]%
{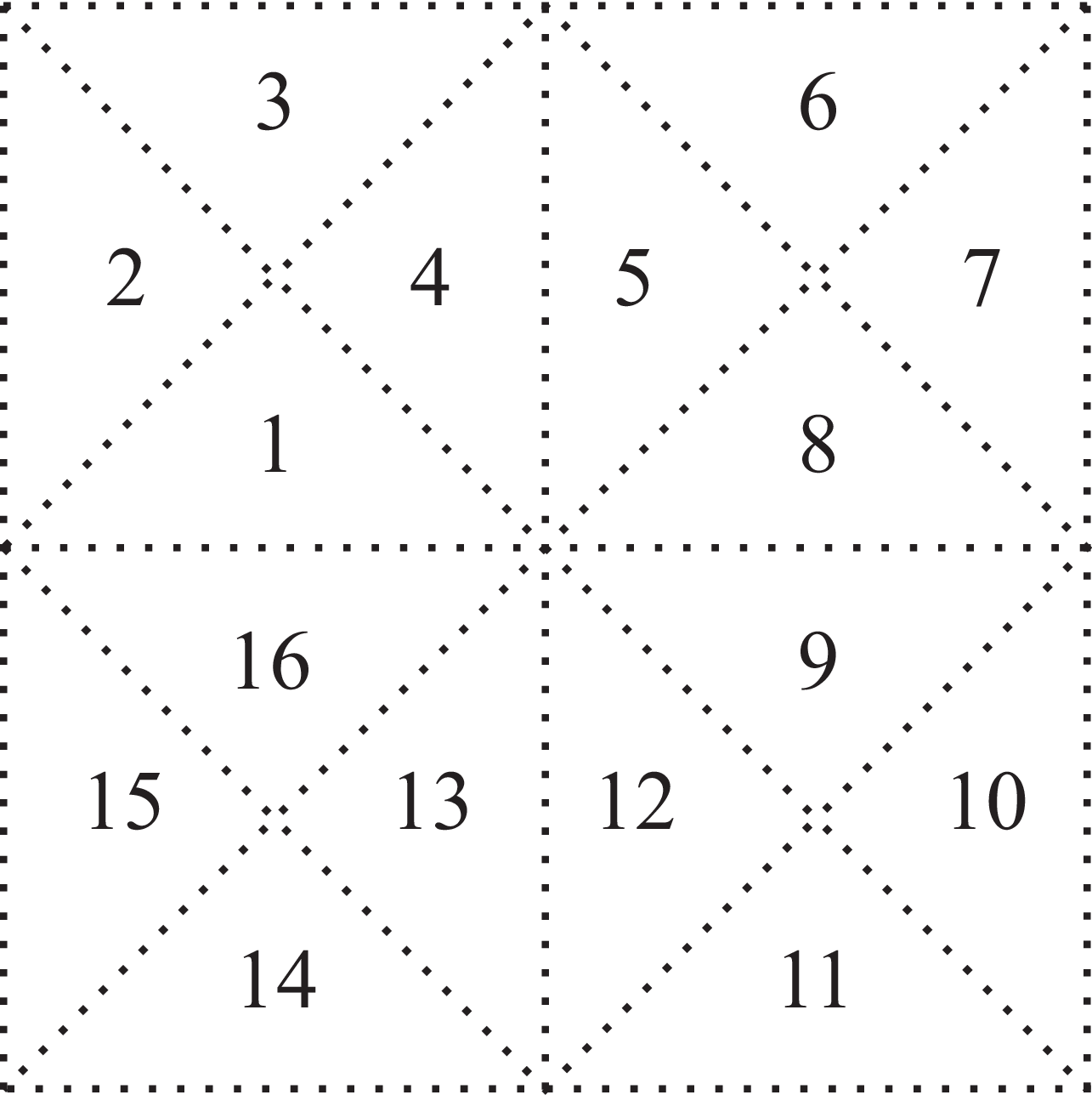}%
\caption{Figure 3. A block of sixteen isosceles right triangular stamps.}%
\end{figure}

A sheet of postage stamps is suitable for stamp folding puzzles if the stamps
are configured in such a way that the sheet folds neatly into a packet of
single stamps. Such a packet unfolds into an edge tessellation in which the
perforations form lines of symmetry. On page 143 of \textit{\cite{Fr}},
Frederickson poses the following conjecture:\textit{ }\emph{Although
triangular stamps have come in a variety of different triangular shapes, only
three shapes seem suitable for [stamp] folding puzzles: equilateral, isosceles
right triangles, and }$60^{\circ}$\emph{-right triangles.} Theorem \ref{one}
confirms Fredickson's Conjecture; indeed, the four suitable edge tessellations
are pictured in Figure 1. We invite the reader to reproduce and fold each of
them into a packet of single stamps. Our folding algorithms appear at the end
of this article.

The question posed at the outset of this article and answered by Theorem
\ref{one}, was motivated by the \textquotedblleft unfolding
technique\textquotedblright\ applied by A. Baxter and the second author to
find, classify and count classes of periodic orbits of a billiard ball in
motion on an equilateral triangular billiard table (see \cite{Ba-Um} for
details). Periodic orbits on polygonal billiard tables of the eight polygonal
types in Theorem \ref{one} unfold as straight line segments in an edge
tessellation. During an REU in 2001, Andrew Baster and students Ethan McCarthy
and Jonathan Eskreis-Winkler applied this unfolding technique to find,
classify, and count classes of periodic orbits on square, rectangular, and
isosceles right triangular billiard tables (see \cite{Ba-Es-Mc}). Presumably,
this technique also applies on billiard tables of the five remaining polygonal types.

Edge tessellations are \emph{wallpaper patterns, }which are tessellations of
the plane with translational symmetries of minimal length in two independent
directions (the group of translational symmetries is discrete). A point $C$ in
a wallpaper pattern is an $n${\textit{-}}\emph{center}{\textit{ }}if the group
of rotational symmetries centered at $C$ is generated by a rotation of minimal
positive rotation angle $\phi_{n}=${$360^{\circ}/n$}.

The students' original proof of Theorem \ref{one} applies the powerful
\emph{Crystallographic Restriction Theorem, }which tightly constrains the
order of a group of rotational symmetries:\ \emph{If }$C$\emph{ is an }%
$n$\emph{-center of a wallpaper pattern, then }$n\in${$\left\{
2,3,4,6\right\}  $ (for a proof, see \cite{Ma} for example). The proof of
Theorem \ref{one} presented here is independent of Crystallographic
Restriction and more geometrically revealing. \smallskip}

\textit{Proof. }We begin the proof of Theorem \ref{one} by constructing a set
$S$ containing the measures of the interior angles of a generating polygon
$G$. Let $V$ be a vertex of $G,$ and let $\theta$ be the measure of the
interior angle at $V;$ then $\theta<180^{\circ}.$ Let $G^{\prime}$ be the
image of $G$ when reflected in an edge of $G$ containing $V.$ Then the
interior angle of $G^{\prime}$ at $V$ has measure $\theta$, and inductively,
the interior angle at $V$ of every copy of $G$ with vertex $V$ has measure
$\theta$ (see Figure 4).Since successively reflecting in the edges of $G$ that
meet at $V$ is a rotational symmetry of angle $2\theta$, the vertex $V$ is an
$n$-center for some $n.$ If $G^{\prime}$ is the rotational image of $G$, then
$\phi_{n}=\theta;$ otherwise $\phi_{n}=2\theta.$ In either case,
$n\theta=360^{\circ}$ for some $n\in\mathbb{N},$ and it follows that every
interior angle of $G$ lies in the set
\[%
\begin{array}
[c]{cl}%
S & =\left\{  x\leq120^{\circ}\mid nx=360^{\circ},\text{ }n\in\mathbb{N}%
\right\}  \smallskip\\
& =\left\{  120^{\circ},90^{\circ},72^{\circ},60^{\circ},51\frac{3}{7}^{\circ
},45^{\circ},40^{\circ},36^{\circ},\ldots,18^{\circ},\ldots\right\}  .
\end{array}
\]
\begin{center}
\includegraphics[
height=1.2047in,
width=1.7071in
]%
{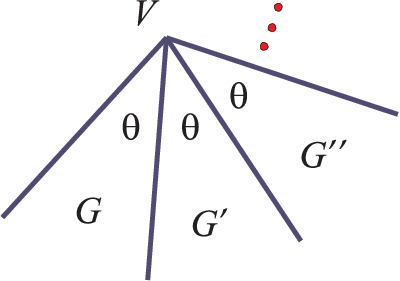}%
\\
Figure 4. Congruent interior angles at vertex $V$ shared by $G$ and images
$G^{\prime}$ and $G^{\prime\prime}.$%
\end{center}

Now suppose that $\theta=120^{\circ};$ then three copies of $G$ share the
vertex $V.$ Let $e$ and $e^{\prime}$ be the edges of $G$ that meet at $V,$ and
labeled so that the angle from edge $e$ to edge $e^{\prime}$ measures
$120^{\circ}$ (see Figure 5). Let $e^{\prime}$ and $e^{\prime\prime}$ be their
respective images under a $120^{\circ}$ rotation. Then $e^{\prime\prime}$ lies
on the bisector of $\angle V$ and is the reflection of $e^{\prime}$ in $e.$ By
a similar argument, if an odd number of copies of $G$ share vertex $V,$ the
bisector of $\angle V$ is a line of symmetry.
\begin{center}
\includegraphics[
height=1.6086in,
width=1.7469in
]%
{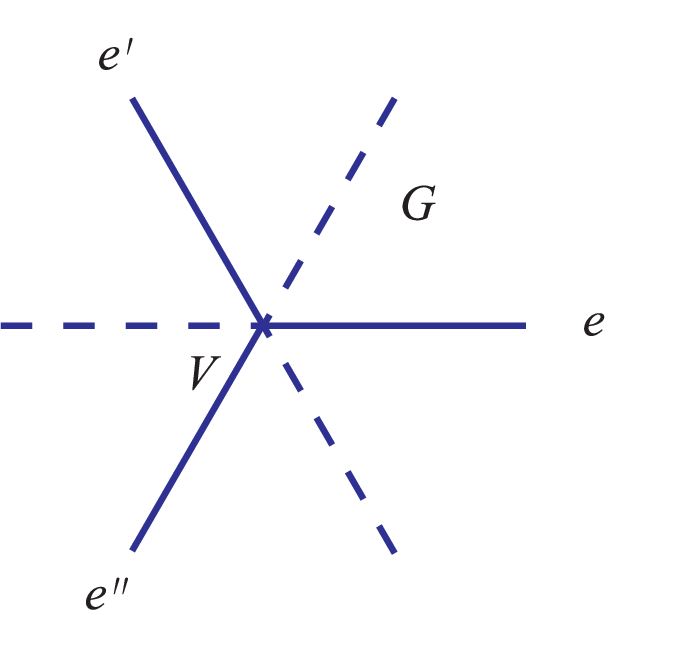}%
\\
Figure 5. A line of symmetry bisects a $120^{\circ}$ interior angle.
\end{center}

Let $g$ be the number of edges of $G.$ Then the interior angle sum
$180^{\circ}\left(  g-2\right)  \leq120^{\circ}g$ implies $g\leq6.$ If $G$ is
a hexagon, it is equiangular since its interior angles measure at most
$120^{\circ}$ and its interior angle sum is $720^{\circ}=6\left(  120^{\circ
}\right)  $. But $G$ is symmetric with respect to each of its interior angle
bisectors by the remark above. Therefore $G$ is a \emph{regular hexagon}.

We claim that $G$ is not a pentagon. On the contrary, suppose $G$ is a
pentagon. Then some interior angle measures $120^{\circ}$ since the interior
angle sum of $540^{\circ}>5\left(  90^{\circ}\right)  $. Choose an interior
angle of $120^{\circ}$ and label the vertex at this interior angle $V.$ Then
$G$ is symmetric with respect to the angle bisector at $V$ and the other
interior angles of $G$ pair off congruently--two with measure $x,$ two with
measure $y.$ Note that the interior angles in one of these pairs are adjacent
(see Figure 6). If $x=y,$ then $x=105^{\circ}\notin S;$ hence $x\neq y.$ If
$x<y,$ then $y>105^{\circ};\ $hence $y=120^{\circ}$ since $y\in S.$ But if
$y=120^{\circ},$ lines of symmetry bisect three interior angles of $G$, in
which case $x=y$ by the adjacency noted above, and $G$ is equiangular with an
interior angle sum of $600^{\circ},$ which is a contradiction.$\bigskip$
\begin{center}
\includegraphics[
height=1.5939in,
width=1.7746in
]%
{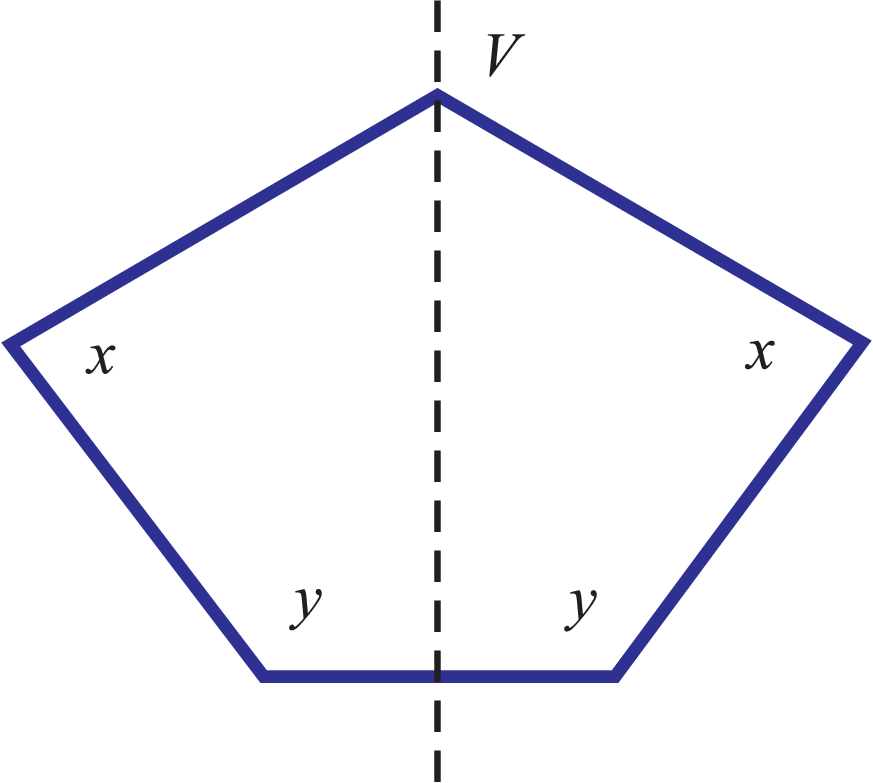}%
\\
Figure 6. A pentagon $G$ with an interior angle of $120^{\circ}$ at $V$.
\end{center}

Now if $V$ is a vertex of $G$, let $m\angle V$ denote the measure of the
interior angle at $V$. Suppose $G$ is a quadrilateral. If $G$ has an interior
angle of $120^{\circ},$ label the vertices $A,B,C,D$ in succession with
$m\angle A=120^{\circ}.$ Then the bisector $s$ of $\angle A$ is a line of
symmetry, $C$ is on $s,$ and $\angle B\cong\angle D$ (see Figure 7). Let
$2x=m\angle C$ and $y=m\angle B,$ and note that $m\angle BAC=60^{\circ}.$ Then
$x\leq60^{\circ}\leq120^{\circ}-x=y.$ Hence the only solutions of
$x+y=120^{\circ}$ with $x\leq y$ and $\left(  x,y\right)  \in S\times S$ are
$\left\{  \left(  30^{\circ},90^{\circ}\right)  ,\left(  60^{\circ},60^{\circ
}\right)  \right\}  .$ Therefore $G$ is either a $120$\emph{-rhombus }or a
$60$-$90$-$120$ \emph{kite. }On the other hand, if the interior angles of $G$
measure at most $90^{\circ},$ then $G$ is equiangular since its interior angle
sum is $360^{\circ},$ and $G$ is a \emph{rectangle}.%

\begin{center}
\includegraphics[
height=1.5273in,
width=1.6276in
]%
{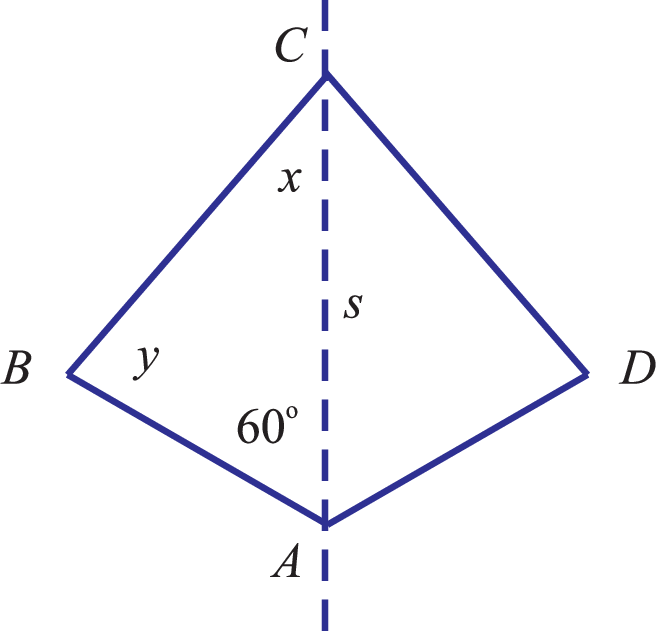}%
\\
Figure 7. A quadrilateral $G$ with an obtuse interior angle.
\end{center}

Finally, suppose $G$ is a triangle. If $G$ has an interior angle of
$120^{\circ}$, then $G$ is a $120$\emph{-isosceles triangle} by symmetry.
Otherwise, let $G=\Delta ABC;$ let $x=m\angle A,$ $y=m\angle B,$ and
$z=m\angle C$.

If $G$ is a right triangle with $z=90^{\circ}$ and $x\leq y,$ then $x,y\in S$
implies $18^{\circ}\leq x\leq y\leq72^{\circ}.$ Hence the only solutions of
$x+y=90^{\circ}$ with $x\leq y$ and $\left(  x,y\right)  \in S\times S$ are
$\left\{  \left(  18^{\circ},72^{\circ}\right)  ,\left(  30^{\circ},60^{\circ
}\right)  ,\left(  45^{\circ},45^{\circ}\right)  \right\}  .$ Furthermore, if
$y=72^{\circ},$ five copies of $G$ share vertex $B$ and the bisector of
$\angle B$ is a line of symmetry, in which case $x=y=90^{\circ},$ which is a
contradiction. Therefore $G$ is either a $60$-\emph{right} or an
\emph{isosceles-right} \emph{triangle}.

If $G$ is an acute triangle with $x\leq y\leq z\leq72^{\circ},$ then
$x=180^{\circ}-\left(  y+z\right)  \geq180^{\circ}-2(72^{\circ})=36^{\circ};$
on the other hand, $x=180^{\circ}-\left(  y+z\right)  \leq180^{\circ
}-2(60^{\circ})=60^{\circ}.$ But $36^{\circ}\leq x\leq60^{\circ}$ implies
$120^{\circ}\leq y+z\leq144^{\circ}.$ Thus if $y\leq z$ and $\left(
y,z\right)  \in S\times S,$ then $\left(  y,z\right)  \in\left\{  \left(
60^{\circ},60^{\circ}\right)  ,\left(  60^{\circ},72^{\circ}\right)  ,\left(
72^{\circ},72^{\circ}\right)  \right\}  $ so that the only solutions of
$x+y+z=180^{\circ}$ with $x\leq y\leq z$ and $\left(  x,y,z\right)  \in
S\times S\times S$ are $\left\{  \left(  36^{\circ},72^{\circ},72^{\circ
}\right)  ,\left(  60^{\circ},60^{\circ},60^{\circ}\right)  \right\}  .$ But
interior angles of $72^{\circ}$ are bisected by lines of symmetry, so the
solution $\left(  36^{\circ},72^{\circ},72^{\circ}\right)  $ is extraneous and
$G$ is an \emph{equilateral triangle. }

This completes the proof of Theorem \ref{one}.\smallskip

We remark that edge tessellations represent 3 of the 17 symmetry types of
wallpaper patterns. Using the labeling defined in \cite{Ma}, general
(non-square) rectangles generate patterns of type $pmm;$ isosceles right
triangles and squares generate patterns of type $p4m$; and the other six
polygons in Theorem \ref{one} generate patterns of type $p6m.$

Here are some explicit algorithms for folding the sheets of stamps in Figure 1
into packets of single stamps. Assume that $T$ is generated by one of the four
non-obtuse polygons identified in Theorem \ref{one}. Choose an infinite strip
$S$ of minimal width bounded by parallel perforation lines $l$ and $m,$ and
\textquotedblleft accordion-fold\textquotedblright\ $T$ onto $S,$ i.e., fold
along $l$ then along $m$ so that $S$ has four leaves configured as a
\textquotedblleft w\textquotedblright, then fold again along $l$ and again
along $m$ so that $S$ has eight \textquotedblleft zig-zag\textquotedblright
leaves, and continue in this manner indefinitely.

Choose a stamp $P$ in $S.$ If $P$ is a rectangle, accordion-fold $S$ onto $P.$
If $P$ is a right triangle, $P$ together with some subset of its images
tessellate a rectangle $R$ of minimal area contained in $S,$ so accordion-fold
$S$ onto $R,$ then fold $R$ onto $P.$ If $P$ is an equilateral triangle, two
of its edges lie in the interior of $S.$ Label these edges $a$ and $b;$ then
its third edge $c$ is contained in the boundary of $S.$ Let $S_{1}$ and
$S_{2}$ be the subsets of $S-P$ bounded by $l,$ $m,$ and $a,$ and by $l,$ $m,$
and $b,$ respectively. Fold $S_{1}$ along $a,$ then along $c,$ then along $b,$
and continue in this manner indefinitely to form an infinite \textquotedblleft
spiral\textquotedblright; similarly, fold $S_{2}$ along $b,$ then along $c,$
then along $a,$ and continue indefinitely.

On the other hand, if $T$ is generated by an obtuse polygon $G,$ it has a
$3$-center $C$ shared by three copies of $G$. Since the interior angle of $G$
at $C$ is bisected by a line of symmetry $l,$ which contains an edge of some
copy of $G,$ folding along $l$ creases the stamp $G.$ Thus $T$ is not suitable
for stamp folding puzzles, and we have established Frederickson's conjecture:

\begin{theorem}
The edge tessellations suitable for stamp folding puzzles are generated by the
four non-obtuse polygons indicated in Theorem \ref{one}.
\end{theorem}

To summarize, we have proved that a polygon generating an edge tessellation is
one of the following eight types: a rectangle; an equilateral, $60$-right,
isosceles right, or $120$-isosceles triangle; a $120$-rhombus; a $60$%
-$90$-$120$ kite; or a regular hexagon. Of these, the four non-obtuse polygons
generate tessellations suitable for stamp folding puzzles; this establishes
Frederickson's Conjecture. Our proof of Frederickson's Conjecture exhibits
explicit algorithms for folding the sheets of stamps in Figure 1 into packets
of single stamps.\medskip

\textbf{\noindent Acknowledgements. }We wish to thank Andrew Hall and Joshua
York for enthusiastically sharing their creative ideas in the early stages of
this project, Natalie Frank for sharing her thoughts on stamp folding
algorithms, and Andrew Baxter, Deirdre Smeltzer, Jim Stasheff and Doris
Schattschneider for offering many helpful editorial suggestions.

\end{document}